\documentclass[12pt]{article}
\usepackage{amsmath,amsfonts,latexsym}

\usepackage{amsfonts}

\def\C{{\Bbb C}}

\newcommand{\qed}{\nopagebreak\par\hspace*{\fill}$\Box$\par\vskip .1mm}
\begin{document}
                                  
\newtheorem{thm}{Theorem}[section]   
\newtheorem{cor}[thm]{Corollary}     
\newtheorem{lem}[thm]{Lemma}
\newtheorem{prop}[thm]{Proposition}
\newtheorem{defn}[thm]{Definition}
\newtheorem{rem}[thm]{Remark}
\newtheorem{ex}[thm]{Example}

\begin{center}
{\LARGE{\bf Classification of irreducible integrable highest weight  modules for
current  Kac-Moody Algebras}}  \\ [5mm] 
{\bf  S. ESWARA  RAO AND   PUNITA BATRA}
\end{center}

\vskip 5mm

{\bf Abstract.} This paper classifies irreducible, integrable highest weight 
modules for
  ``current Kac-Moody Algebras'' with finite dimensional weight
  spaces. We prove that these modules turn out to be  modules of
  appropriate direct sums of finitely many copies of  Kac-Moody Lie
  algebras.



\section{Introduction}
\label{1.}
Weight modules of current algebras  have
 been studied
 by several authors recently. A complete classification of the simple weight modules for current algebras has been given by
 D. Britten, M. Lau and L. Frank ([4]). In ([2]), M. Bennett, A. Berenstein, V. Chari, A. Khoroshkin and S. Loktev studied
 Macdonald polynomials and BGG reciprocity for current algebras.
The purpose  of this paper is to classify irreducible, integrable,  highest 
weight modules for current Kac-Moody algebras $\mathfrak g_{A}$ with finite-dimensional weight
  spaces . Let   $\mathfrak g$ be a Kac-Moody Lie algebra.
Let $\mathfrak h$ be a Cartan subalgebra of 
$\mathfrak g$. Let   $\mathfrak g^{\prime} = [\mathfrak g, \mathfrak g]$ be the derived algebra of $\mathfrak g$. Let $\mathfrak h^{\prime} = \mathfrak h \cap \mathfrak g^{\prime}$ and $\mathfrak h^{\prime \prime}$ be a complimentary subspace of $\mathfrak h^{\prime}$ in $\mathfrak h$. Then $\mathfrak h = \mathfrak h^{\prime} \oplus 
\mathfrak h^{\prime \prime}$. 
Let $A$ be a commutative, associative, unital algebra which is finitely generated over $\C$.
Let $\mathfrak g_{A} = (\mathfrak g^{\prime}  \otimes A) \oplus \mathfrak h^{\prime \prime}  $. These algebras include the finite-dimensional simple Lie algebras, the multi-loop algebras $\mathfrak g \otimes \C[t_1^{\pm 1} , t_2^{\pm 1},..., t_n^{\pm 1}]$, where $\mathfrak g$ is a finite dimensional simple Lie algebra over $\C$ and current algebras (Example 2.4 in ([6])). 
The most interesting modules of these algebras are the integrable modules(where the real root spaces acts locally nilpotently). 
 Irreducible, integrable modules for multi-loop algebras with finite dimensional weight spaces have been classified by S. Eswara Rao ([10]). Pal and Batra ([7])
have studied the representations of graded multiloop algebras.
Representations of twisted multiloop Lie algebras have been studied by Batra ([1]). 

In the Laurent polynomial ring, if there is only one variable i. e. when $n = 1$, then this multi-loop Lie algebra is called a loop algebra. Representations of loop algebras have been given by Rao ([8]).
The universal central extension of multi-loop algebra $\mathfrak g \otimes \C[t_1^{\pm 1} , t_2^{\pm 1},..., t_n^{\pm 1}]$ is called a toroidal Lie algebra ([5]), 
([9]). Berman and Billig ([3]) have constructed integrable, irreducible  modules for a more general toroidal Lie algebra.
Classification of
  irreducible integrable modules for toroidal Lie algebras
with  finite dimensional weight spaces has been done by Rao ([11]).

In Section 2, we define current Kac-Moody algebra $\mathfrak g_{A}$ and its
highest weight  module(Definition 2.1). Using $\psi$ defined in Definition 2.1, we get the unique 
irreducible $\mathfrak g_{A}$-module $V(\psi)$. Our main problem is to find out
when $V(\psi)$ has finite-dimensional weight spaces. In Lemma 2.3, it has been shown that 
 $V(\psi)$ has finite dimensional weight spaces with respect to subalgebra $\mathfrak h$ if and only if there exists a cofinite ideal $I$ of $A$ such that $\psi(\mathfrak h^{\prime} \otimes I) = 0$.
Such an ideal $I$ of $A$ exists if and only if $I$ satisfies
$\mathfrak g^{\prime} \otimes I.  V(\psi)= 0$(Proposition 2.4), where 
 $\mathfrak g^{\prime}$ is the derived algebra of $\mathfrak g$. 

In Section 3, we define radical of an ideal $I$ in a commutative ring and also
define integrable modules for  $\mathfrak g_{A}$. In the main Theorem 3.3, we show if $V(\psi)$ is an integrable  $\mathfrak g_{A}$-module with finite dimensional weight spaces, then 
$V(\psi)$ is a module for 
 $\oplus_{k-times}\mathfrak g^{\prime} \oplus {\mathfrak h}^{\prime \prime}$, which is a direct sum of $k$-copies of Kac-Moody Lie algebras. It is well known that 
 $V(\psi) \cong V(\lambda_1) \otimes V(\lambda_2) \otimes.... \otimes  V(\lambda_k)$, where $\lambda_1, \lambda_2,... ,\lambda_k$ are dominant integral weights of $\mathfrak g$.

\section{Preliminaries}
\label{2.}
 Let  $\mathfrak g$ be any Kac-Moody Lie algebra. Let $\mathfrak h$ be a Cartan subalgebra of 
$\mathfrak g$. Let   $\mathfrak g^{\prime} = [\mathfrak g, \mathfrak g]$ be the derived algebra of $\mathfrak g$. Let $\mathfrak h^{\prime} = \mathfrak h \cap \mathfrak g^{\prime}$ and $\mathfrak h^{\prime \prime}$ be a complimentary subspace of $\mathfrak h^{\prime}$ in $\mathfrak h$. Then $\mathfrak h = \mathfrak h^{\prime} \oplus 
\mathfrak h^{\prime \prime}$. 
Let $A$ be a commutative, associative, unital algebra which is finitely generated over $\C$. Let   $\mathfrak g^{\prime}=  
N^{-} \oplus \mathfrak h^{\prime}  \oplus N^{+}$ be the standard decomposition into negative root spaces, positive root spaces and a $\mathfrak h^{\prime}$.  Let $\mathfrak g_{A} = \mathfrak g^{\prime} \otimes A \oplus \mathfrak h^{\prime \prime}$.
Then $\mathfrak g_{A} = ({N^{-} \otimes A}) \oplus ({ \mathfrak h^{\prime} \otimes A) \oplus \mathfrak h^{\prime \prime}  }    \oplus ({ N^{+} \otimes A}) $ and it
can be made into a Lie algebra by defining 
$$[ X \otimes a,  Y \otimes b]= [ X,Y]  \otimes ab,$$
$$[ X \otimes a, h] = [X,h]  \otimes a,$$
for all $X, Y \in \mathfrak g^{\prime}, h \in \mathfrak h^{\prime \prime}$  and $a, b \in A$.
Then $\mathfrak g_{A}$ is called a current Kac-Moody Lie algebra.

\begin{defn}
\label{2.1.}
A  ${\mathfrak g}_{A}$-module $W$ is said to be a highest weight  module if there exists a weight vector $v$ in $W$ such that

(1) $U({\mathfrak g}_{A}) v = W$,

(2)  ${ N^{+} \otimes A}. v = 0$,

(3) there exists a linear  map $\psi : ({ \mathfrak h^{\prime} \otimes A}) \oplus
\mathfrak h^{\prime \prime} \rightarrow \C$  such that \\
    $h.v =\psi(h) v $ for all $h \in ({ \mathfrak h^{\prime} \otimes A}) \oplus
\mathfrak h^{\prime \prime}$.
\end{defn}

Let $\psi$ be as in Def 2.1 and let $\C(\psi)$ be a one-dimensional vector space. Let ${N^{+} \otimes A}$ act trivially on $\C(\psi)$ and let $({ \mathfrak h^{\prime} \otimes A}) \oplus
\mathfrak h^{\prime \prime}$ act by $\psi$. Now consider the induced module for 
 ${\mathfrak g}_{A}$ given by
$$ M(\psi) = U({\mathfrak g}_{A}) \otimes_{B} \C(\psi),$$ where $B = ({ N^{+} \otimes A}) \oplus ({ \mathfrak h^{\prime} \otimes A}) \oplus
\mathfrak h^{\prime \prime}.$

By standard arguments one shows that  $M(\psi)$ is an $\mathfrak h$-weight module  and has a unique irreducible quotient, which we will denote by $V(\psi)$. 
Our problem is to find out when $V(\psi)$ has finite dimensional weight spaces.
 By definition an ideal $I$ in $A$ is cofinite if and only if  $A/I$ is a finite dimensional vector space over $\C$.

\begin{lem}
\label{2.2.}
Suppose $I$ and $J$ are cofinite ideals of $A$, Then $IJ$ is also a cofinite ideal of $A$. In particular, if $I_1, I_2,...., I_n$ are cofinite ideals of $A$, then $I_1I_2....I_n$ is also a cofinite ideal in $A$.
\end{lem}
{\bf Proof.}

Since  $A$ is a finitely generated algebra over the  field $\C$, we note that $A$ is always Noetherian.
Since $I$ is a cofinite ideal of $A$, so by definition
the quotient ring $A/I$ is a finite dimensional vector space over $\C$.
The quotient ring $A/I$ is a finite dimensional vector space over $\C$
 if and only if it is Artinian, i. e. it is 
Noetherian and of Krull dimension 0.
 Krull dimension 0 means by definition every prime ideal is maximal.
 Since  $I$ and $J$ are ideals of  $A$ and $A$ is Noetherian, hence  $A/IJ$ is
also  Noetherian(being the homomorphic image of $A$). Therefore to prove that $A/IJ$ is
again cofinite, it is enough to prove that every prime ideal in $A/IJ$ is maximal. Now, every
prime ideal  $P$ in $A/IJ$ corresponds to  a prime ideal in $A$ which contains the product ideal $IJ$
and hence it must contain one of the ideal $I$ or $J$. Therefore this prime ideal corresponds to a
prime ideal either in $A/I$ or in $A/J$. This proves that it must be maximal 
by assumption of 
confinite on $I$ and $J$. This completes the proof.
{\qed}

\begin{lem}
\label{2.3.}
Let $\psi$ be as in Definition 2.1. Then the irreducible  ${\mathfrak g}_{A}$-module $V(\psi)$ has finite dimensional weight spaces with respect to $\mathfrak h$ if and only if there exists a cofinite ideal $I$ of $A$ such that $\psi(\mathfrak h^{\prime} \otimes I) = 0$.
\end{lem}

First we prove the following Proposition.

\begin{prop}
\label{2.4.}
There exists an ideal $I$ of $A$ such that $\psi(\mathfrak h^{\prime} \otimes I) = 0$ if and only if $I$ satisfies
$\mathfrak g^{\prime} \otimes I.  V(\psi)= 0$.
\end{prop}
{\bf Proof.} 
Suppose $\psi(\mathfrak h^{\prime} \otimes I) = 0$  for some  ideal $I$ of $A$.  Let 
$\Delta^+$ be a positive root system for $\mathfrak g$. Let ${\alpha}_i$ be a simple root in $\Delta^+$. Let$X_{{\alpha}_i}$, $X_{- {\alpha}_i}$ be  root vectors for the roots ${\alpha}_i$ and
$- {\alpha}_i$ such that $[X_{{\alpha}_i}, X_{- {\alpha}_i}] = h_i$. Consider, for a highest weight vector $v$, the set $(X_{- {\alpha}_i} \otimes I) v$, which is contained in the weight space $V_{\psi - {\alpha}_i}(\psi)$

{\bf Claim 1}$$(X_{- {\alpha}_i }\otimes I) v = 0.$$

In order to prove Claim 1 we will prove that
$$(X_{\alpha} \otimes A)(X_{- {\alpha}_i }\otimes I) v = 0 \ for \ all   \ {\alpha}\  positive \ roots.  \ \ \ \ (*)$$  
It is sufficient to prove  $(*)$ for simple roots. We note that if $j \neq i$, then $$(X_{{\alpha}_j} \otimes A)(X_{- {\alpha}_i }\otimes I) v = (X_{- {\alpha}_i }\otimes I)(X_{{\alpha}_j} \otimes A)v= 0,$$
since $v$ is a highest weight vector.
If $j =  i$, then $$(X_{{\alpha}_j} \otimes A)(X_{- {\alpha}_i }\otimes I) v = ( h_i \otimes I) v = 0,$$ by hypothesis. Hence $(X_{- {\alpha}_i} \otimes I) v$
is a highest weight space.
So  $(X_{- {\alpha}_i} \otimes I)v$ generates a proper submodule of the irrducible module $V(\psi)$. Hence $(X_{- {\alpha}_i }\otimes I) v = 0.$

{\bf Claim 2}$(X_{- {\alpha} }\otimes I) v = 0 \ for \ all  \ \ {\alpha} \in \Delta^+ .$
The claim is true for any simple root by claim 1. This reduces the proof to a proof by induction on the height of $\alpha$. The claim is true for  $\alpha$ such that height $\alpha = 1$. Let $\alpha_i$ be a simple root in $\Delta^+$, then 
$$(X_{{\alpha}_i} \otimes A)(X_{- {\alpha} }\otimes I) v = [X_{{\alpha}_i}, X_{- {\alpha} }] \otimes I v + (X_{- {\alpha} }\otimes I) (X_{{\alpha}_i} \otimes A) v.$$
The second term is zero as $v$ is a highest weight vector. The first term is zero by induction. Since $(X_{- {\alpha} }\otimes I) v$ is killed by $X_{{\alpha}_i} \otimes A$ for  $\alpha_i$ simple, it is easy to see that 
 $(X_{- {\alpha} }\otimes I) v$ is a highest weight space of weight $
\psi -  {\alpha}$. But $V(\psi)$ is an irreducible highest weight module and hence $(X_{- {\alpha} }\otimes I) v = 0$. 
So  $$(X_{- {\alpha} }\otimes I) v = 0    \ for \ all \ {\alpha} \in \Delta^+ ,$$  $$({\mathfrak h}^{\prime} \otimes I) v = 0 \  by \  hypothesis,$$  
  $$(X_ {\alpha} \otimes I) v = 0   \ for \ all \  {\alpha} \in \Delta^+ ,\  by \ definition.$$

So we have proved   ${\mathfrak g}^{\prime} \otimes I v = 0$. Now consider
$$W = \{ w \in V(\psi) |  {\mathfrak g}^{\prime} \otimes I w = 0 \}.$$
It is easy to see that $W$ is a ${\mathfrak g}_A$-submodule of $V(\psi)$.
Since $v  \in W$, so $W$ is nonempty. So $W$ is a nonzero submodule of
 $V(\psi)$.  But $V(\psi)$ is irreducible. Hence $W = V(\psi)$. This implies
$\mathfrak g^{\prime} \otimes I.  V(\psi)= 0$.

Conversely if $\mathfrak g^{\prime} \otimes I.  V(\psi)= 0$ for some $I$.
We know that ${\mathfrak h}^{\prime} \otimes I \subset  {\mathfrak g}^{\prime} \otimes I$ and $\C v \subset  V(\psi)$. Hence ${\mathfrak h}^{\prime} \otimes I v = 0$.
 
{\qed}
{\bf Proof of Lemma 2.3.} 
Suppose $I$ is a cofinite ideal of $A$ such that $\psi({\mathfrak h}^{\prime} \otimes I) = 0$. Then by Proposition 2.4, 
$$(\mathfrak g^{\prime} \otimes A) v = (\mathfrak g^{\prime} \otimes A/I) v = U(N^{-} \otimes A/I) v.$$
So $V(\psi)$ is a module for $(\mathfrak g^{\prime} \otimes A/I)$ and hence 
 $V(\psi)$ has finite dimensional weight spaces using PBW theorem.

Conversely suppose $V(\psi)$ has finite dimensional weight spaces. Let $\alpha_i$ 
be a simple root in  $\Delta^+$ and $X_{{\alpha}_i}$ be a root vector for the root
${{\alpha}_i}$. Consider $X_{-{\alpha}_i}$ be a root vector for the root
${-{\alpha}_i}$ such that $[X_{{\alpha}_i}, X_{-{\alpha}_i}] = h_i$. Consider
$$I_i = \{ a \in A | X_{-{\alpha}_i}  \otimes a v = 0 \}.$$
{\bf Claim} $I_i$ is a cofinite ideal of $A$. \\
First we show that $I_i$ is an  ideal of $A$.
We have for $a \in I_i$ and  $b \in A$
  $$
 \begin{array}{lll}
0  &=& ( h_i \otimes b) ( X_{-{\alpha}_i}  \otimes a ) v \\
      &=&   ( X_{-{\alpha}_i}  \otimes a )  ( h_i \otimes b)v + [ h_i ,  X_{-{\alpha}_i}]  \otimes ab v \\
      &=&   \psi( h_i \otimes b) ( X_{-{\alpha}_i}  \otimes a ) v - 2X_{-{\alpha}_i}  \otimes ab v \\
       &=&  -2 X_{-{\alpha}_i}  \otimes ab v.$$
\end{array}
$$
This shows that  $ab \in I_i$, so $I_i$ is an  ideal of $A$. Now we observe that
elements of the form $X_{-{\alpha}_i}  \otimes a v, a \in A$ are contained in the weight space 
$V_{\psi - {\alpha}_i}(\psi)$, which is a finite dimensional space, by hypothesis.
Let dim$V_{\psi - {\alpha}_i}(\psi) = k$
{\bf Claim}  dim$A/I \leq k$. \\
Let $a_1, a_2,... ,a_n \in A$, where $n > k$.
Then   $$X_{-{\alpha}_i}  \otimes a_1 v, X_{-{\alpha}_i}  \otimes a_2 v,... , X_{-{\alpha}_i}  \otimes a_n  v \in V_{\psi - {\alpha}_i}(\psi).$$
So they must be linearly dependent vectors. So  for $m_1, m_2, \cdots, m_n \in \C$, we have
$$X_{-{\alpha}_i}  \otimes \sum_{i=1}^{n} m_ia_i v = 0.$$
Hence $$ \sum_{i=1}^{n} m_ia_i  \in I_i .$$
This shows that dim$A/I_i \leq k$. So  $I_i$ is a non-zero cofinite ideal of $A$. Let rank of ${\mathfrak g} = l$ and consider $I =  I_1^{2} I_2^{2}....I_l^{2}$, which is also cofinite by Lemma 2.2.
So we have $ X_{-{\alpha}_i}  \otimes I_i v = 0$  for  all  simple $i$,  
 and $ X_{{\alpha}_i}  \otimes I_i v = 0$ by definition of highest weight module.
So it follows  that $ h_i  \otimes I_i^{2} v = 0$  for  all   $i$.
Since $I =  I_1^{2} I_2^{2}....I_l^{2} \subset I_i^{2}$ for all $1 \leq i \leq l$, so  $ h_i  \otimes I v = 0$  for  all   $i$  implies that  ${\mathfrak  h^{\prime}}  \otimes I v = 0$.
{\qed}

\section{Integrable ${\mathfrak g}_{A}$-Modules }

\begin{defn}
\label{3.1.}
The radical of an ideal $I$ in a commutative ring $A$, denoted by $\sqrt I$ is defined as

$$\sqrt I = \{ a \in A | a^n \in I \ for \ some \ positive \ integer \ n \}.$$
It is easy to see that $\sqrt I $ is an ideal itself containing $I$.
\end{defn}

\begin{defn}
\label{3.2.}
A module $V$ of ${\mathfrak g}_{A}$ is called integrable if

1) $V = \bigoplus_{\lambda \in {\mathfrak
   h}^* } V_{\lambda},$  where     
  $$V_{\lambda} = \{ v \in V  |  \  h v = \lambda(h) v  \ for \ all \
 h \in {\mathfrak
   h} \},$$

2) For all real roots $\alpha$ of ${\mathfrak g}$ and for all $a \in A$ \
 and \ \
  $v \in
   V$, there exists an integer $k= k(\alpha, v)$ such that
   $(X_{\alpha} \otimes a)^k v = 0$, where
   $X_{\alpha}$ is the root vector corresponding to root $\alpha$.
\end{defn}

\begin{thm}
\label{3.3.}
Let   $V(\psi)$ be an integrable  module for ${\mathfrak g}_{A}$ with finite-dimensional weight spaces. Then $V(\psi)$ is  a module for $(\mathfrak g^{\prime} \otimes A/J ) \oplus {\mathfrak h}^{\prime \prime}$, where $J$ is a radical ideal of $A$.
In particular $V(\psi)$ is a module for $\oplus_{k-times}\mathfrak g^{\prime} \oplus {\mathfrak h}^{\prime \prime}$.
\end{thm}
{\bf Proof.} 
Since $V(\psi)$ has finite-dimensional weight spaces, so by Lemma 2.3, there exists a 
cofinite ideal $I$ of $A$ such that $\psi(\mathfrak h^{\prime} \otimes I) = 0$.
We define a cofinite ideal $\sqrt I$ as in Defiinition 3.1 such that
$$ I \subset \sqrt I,$$
which is a radical ideal of $A$, call it $J$.

{\bf Claim 1}  $\mathfrak g^{\prime} \otimes J$ is zero on  $V(\psi)$. \\
Consider the map
$$\phi : \mathfrak g^{\prime} \otimes A/I \rightarrow \mathfrak g^{\prime} \otimes A/J.$$
First we observe that Ker$\phi =  \mathfrak g^{\prime} \otimes J/I$ is solvable.
Since $J/I$ is a vector subspace of $ A/I$  and  dim$A/I  < \infty$, so $J/I$ is a finite dimensional vector subspace of $ A/I$. Let $a_1 + I, a_2 + I,...., a_k + I$ be 
a basis for 
$J/I$. Then there exists some positive integer $N$ sufficiently large such that 
$a_i^N \in I  \ \forall  \ i = 1,..., k$. \\
Consider
$$[X_{i_k}  \otimes a_{i_k},....[X_{i_2}  \otimes a_{i_2}, X_{i_1}  \otimes a_{i_1}]..]
= [X_{i_k}..[X_{i_2}, X_{i_1}] ..] \otimes a_{i_k}..a_{i_2}a_{i_1},$$
where $a_{i_1}, a_{i_2},..., a_{i_k} \in \{a_1,...,a_k \}$.
For sufficiently many $a_i$'s the product $a_{i_k}..a_{i_2}a_{i_1}$ will be in $I$. Hence  Ker$\phi =  \mathfrak g^{\prime} \otimes J/I$ is a solvable Lie algebra.
 Let $\alpha_i$ be a simple root of  $\mathfrak g$. Let $X_{{\alpha}_i}, X_{-{\alpha}_i}, 
h_{{\alpha}_i}$ be a $sl_2$-triplet and ${\mathfrak g_{{\alpha}_i}}$ be the span of
 $X_{{\alpha}_i}, X_{-{\alpha}_i}$ and
$h_{{\alpha}_i}$ and consider $\overline {\mathfrak g_{{\alpha}_i}} = {\mathfrak g_{{\alpha}_i}} \otimes  A/I$, which is finite-dimensional. Let $W$ be $\overline {\mathfrak g_{{\alpha}_i}}$-module generated by $v$ i.e. $W = U(\overline {\mathfrak g_{{\alpha}_i}})v$. Since 
 $V(\psi)$ is integrable and $h_{{\alpha}_i} \otimes  A/I$ acts as scalars on $v$, we conclude that $W$ is finite-dimensional. By restricting the action of the solvable Lie algebra $\mathfrak g_{{\alpha}_i} \otimes J/I$ to $W$, we get a vector
$w \in W$(using Lie's Theorem) such that $\mathfrak g_{{\alpha}_i} \otimes J/I$ acts as scalars on $w$. Using Proposition 2.1 ([10]), we get that
 $$\mathfrak g_{{\alpha}_i} \otimes (J/I) w = 0 \ \ \  \ \ \ \  \    \ \ \ \  \    \         (1)$$ 

{\bf Claim 2} $$\mathfrak g_{{\alpha}_i} \otimes (J/I) v = 0.$$ 
Using $W$, we have $Y \in U(X_{-{\alpha}_i} \otimes  A/I)$ such that $Yv= w$. Since  $V(\psi)$ is irreducible, there exists $X \in U({\mathfrak g}_{A})$ such that
$$Xw = v= XYv.$$ Without loss of generality, we can assume $w$ is a weight vector. Here $X= X_{-}hX_{+}, X_{+} \in U({N^{+} \otimes A}), X_{-} \in U({N^{-} \otimes A}), h \in U( \mathfrak h^{\prime} \otimes A)$.
 But by weight argument $X= hX_{+}$, $X_{-}$ has to be a constant.
So $$Xw = h_{i_1}a_{i_1}^{\prime}.... h_{i_r}a_{i_r}^{\prime}X_{{\alpha}_i}a_{1} ...X_{{\alpha}_i}a_{l} w = v, \ \ \ \     \      \          \ (2)$$  where
$h_{i_1},...,h_{i_r}  \in \mathfrak h^{\prime}, a_{i_1}^{\prime}, .... ,a_{i_r}^{\prime}, a_{1}, ..., a_{l} \in A$. 
We see that by (1)$$\mathfrak g_{{\alpha}_i} \otimes J  U( \mathfrak h^{\prime} \otimes A) w = 0. \ \  \ \   \     \ \     \    \ \ \ \  \ \ \ \ \  \  (3)$$ 

Next 
$$\mathfrak g_{{\alpha}_i} \otimes J X_{{\alpha}_i}a_{1}w = X_{{\alpha}_i}a_{1} \mathfrak g_{{\alpha}_i} \otimes J w + [ X_{{\alpha}_i}, \mathfrak g_{{\alpha}_i}]  \otimes a_{1} J w = 0,$$
as both terms are zero by (1). Now by induction on $l$, we get that 
 $\mathfrak g_{{\alpha}_i} \otimes J X_{{\alpha}_i}a_{1}.....X_{{\alpha}_i}a_{l} w = 0.$
From (1) and (3), it follows that 
$\mathfrak g_{{\alpha}_i} \otimes J X w = \mathfrak g_{{\alpha}_i} \otimes J v = 0.$
This proves Claim 2.

Since  $\mathfrak h_{{\alpha}_i} \otimes J v  \subset \mathfrak g_{{\alpha}_i} \otimes J  v$, so $\mathfrak h_{{\alpha}_i} \otimes J v = 0$ for each simple root ${\alpha}_i$. Hence $\mathfrak h^{\prime} \otimes J v = 0.$ 
So by Proposition 2.4,
${\mathfrak g}^{\prime} \otimes J  = 0$ on $V(\psi)$.  This proves Claim 1.

It is a standard result that 
 $J = M_1 \cap M_2 \cap ...\cap M_k $(finite intersection), \\
where each $ M_i$ is a maximal ideal and these are all distinct ideals. By the Chinese
remainder theorem
$$ A/J \cong A/M_1 \oplus A/M_2 \oplus ...\oplus A/M_k. $$
Since $A/M_i$ is a field, so $A/M_i \cong \C, \forall i = 1,..., k$.
So $$( \mathfrak g^{\prime} \otimes A/J ) \oplus {\mathfrak h}^{\prime \prime} \cong \oplus_{k-times} \mathfrak g^{\prime}  \oplus {\mathfrak h}^{\prime \prime}, $$
 which is a direct sum of $k$-copies of Kac-Moody Lie algebras. It is well known that 
 $V(\psi) \cong V(\lambda_1) \otimes V(\lambda_2) \otimes.... \otimes  V(\lambda_k)$, where $\lambda_1, \lambda_2,... ,\lambda_k$ are dominant integral weights of $\mathfrak g$.

{\qed}

{\bf Acknowledgments . }
We would like to thank Professor Dilip Patil for some consultation in Commutative Algebra.

\begin{center}
{\bf References}
\end{center}
\begin{enumerate}

\item[{[1]}]  P. BATRA,  Representations of twisted multiloop  Lie algebra,
J.  Algebra {\bf  272} (2004), 404-416.
\item[{[2]}]   M. BENNET, A. BERENSTEIN,  V. CHARI, A. KHOROSHKIN AND  S. 
LOKTEV,
  Macdonald polynomials and BGG reciprocity for current algebras,  Selecta Math. (N.S.) {\bf 20}  (2014),  no. 2, 585 ???607. 
 \item[{[3]}]    S. BERMAN  AND Y. BILLIG,  Irreducible representations for
toroidal Lie algebras,  J.   Algebra {\bf  221} (1999), 188-231.
 \item[{[4]}]   D. BRITTEN, M. LAU AND  L. FRANK,  Weight modules for current algebras, preprint 2014.
 arXiv:1411.3788v1 
\item[{[5]}]    R. V. MOODY, S. E. RAO AND T. YOKUNUMA,  Toroidal Lie algebra and vertex 
representations,    Geometriae Dedicata  {\bf  35} (1990), 283-307.
\item[{[6]}] E. NEHER, A. SAVAGE AND P. SENESI,  Irreducible finite-dimensional  representations of equivariant algebras, Trans.   Amer.
 Math. Soc.
{\bf 364}  (2012),   2619-2646.
\item[{[7]}] T. PAL AND P. BATRA,  Representations of graded multiloop 
algebras, Commun. Algebra  {\bf  38} (2010), 49-67.
\item[{[8]}]    S. E. RAO,  On representations of loop algebras,
   Commun.  Algebra {\bf 21} (1993), 2131-2153. 
\item[{[9]}]   S. E. RAO AND R. V. MOODY,                  Vertex
representations for $n-$toroidal Lie algebras and a generalization of
the Virosoro algebra, Commun.    Math.  Physics
{\bf 159}
        (1994), 239-264.
\item[{[10]}]    S. E. RAO,   Classification of
  irreducible integrable modules for multi loop algebras with finite 
dimensional weight spaces,  J.  Algebra {\bf 246}  (2001), 215-225.
\item[{[11]}]    S. E. RAO,   Classification of
  irreducible integrable modules for toroidal Lie algebras
with  finite dimensional weight spaces,  J.  Algebra
 {\bf 227 } (2004),
318-348.

\end{enumerate}

S. ESWARA RAO \\
SCHOOL OF  MATHEMATICS \\
TATA INSTITUTE OF FUNDAMENTAL RESEARCH \\
 HOMI BHABHA ROAD \\
MUMBAI 400 005 \\
INDIA \\
senapati@math.tifr.res.in \\

PUNITA BATRA \\
HARISH  CHANDRA RESEARCH INSTITUTE \\
CHHATNAG ROAD \\
JHUNSI \\
ALLAHABAD 211 019 \\
INDIA \\
batra@hri.res.in

\end{document}